\documentclass[11pt]{article}
\usepackage{amsmath}
\usepackage{amssymb}
\usepackage{amsthm}
\usepackage{epsf}

\newcounter{AbcT}

\newtheorem {Theorem}    {Theorem}[section]

\newtheorem {Lemma}      [Theorem]    {Lemma}

\newtheorem {Proposition}[Theorem]    {Proposition}
\newtheorem {Claim}      [Theorem]    {Claim}

\newcommand {\Heads}[1]   {\smallskip\pagebreak[1]\noindent{\bf
#1{\hskip 0.2cm}
}}
\newcommand {\Head}[1]    {\Heads{#1:}}

\newcommand {\proofs}     {\proof}
\newcommand {\prooft}[1]  {\Head{Proof #1}\nopagebreak[2]}

\newcommand  {\QED}    {\def\qedsymbol{$\blacksquare$}\qed}

\newcommand{\Z}{{\bf{Z}}}
\newcommand{\ignore}[1]{}

  \def\E{{\bf{E}}}
\def\P{{\bf{P}}}
\def\V{{\bf{V}}}

\def\R{\hbox{I\kern-.2em\hbox{R}}}

\def\C{{\cal{C}}}

\def\G{{\cal{G}}}

\def\|{\, | \, }
\def\v0{{\bf 0}}

\def\0{\hat{0}}
\def\1{\hat{1}}

\def\path{{\tt path}}

\def\lam{\lambda}

\def\phi{\varphi}

\def\be{\begin{equation}}
\def\ee{\end{equation}}

\def\dir{{\cal E}}

\begin{document}
\title{On the mixing time of simple random walk on the super critical percolation cluster}
\author{Itai Benjamini and  Elchanan Mossel \\ Weizmann Institute and Microsoft Research}
\maketitle

\begin{abstract}

We study the robustness under perturbations of mixing
times, by studying mixing times of random walks in percolation clusters inside boxes in $\Z^d$.
We show that for $d \geq 2$ and $p > p_c(\Z^d)$,
the mixing time of simple random walk on the largest cluster inside
$\{-n,\ldots,n\}^d$ is $\Theta(n^2)$ - thus the mixing time is
robust up to constant factor.
\end{abstract}

\section{Introduction}

An important parameter of random walks on finite graph is the mixing time
of the random walk. We refer the reader to \cite{Al} for background and many references,
or to subsection \ref{subsec:background} for terminology used in this paper.

It is natural to study the robustness of the mixing time
under perturbations.
In this paper we address this issue by studying the effect of {\sl random} perturbations
of the underlying graph on the mixing
times of simple random walk inside boxes in $\Z^d$.

A classical way to perturb the lattice $\Z^d$ is by performing super-critical percolation.
See \cite{g} for introduction and many references,
or subsection \ref{subsec:background} for terminology.

In this paper we study the mixing
time of simple random walk on the
largest super-critical percolation cluster inside $\{-n,\ldots,n\}^d$.
We show that for all $d \geq 2$, the mixing time of the random walk on the perturbed
box, is up to constant, the same as the mixing time of simple random walk on the original box.

Below we use the notation $g(n) = O(f(n))$ to indicate that
$\limsup g(n)/f(n) < \infty$. We will write $f(n) = \Omega(g(n))$
when $g(n) = O(f(n))$, and $f(n) = \Theta(f(n))$ when both
$f(n) = O(g(n))$ and $g(n) = O(f(n))$.

\subsection{Main results}

Let $G=(V,E)$ be a finite graph with no loops.
Consider a simple random walk on $G$. In order to avoid periodicity of the random walk,
we work with the continuous time random walk.
Look at the markov semi-group generated by the matrix $(Q_{x,y})_{x \in V,y \in V}$
where \begin{equation} \label{eq:transitions} Q_{x,y} = \left\{ \begin{array}{ll}
            \frac{1}{\deg(x)} & \mbox{if } x \sim y, \\
            -1 & \mbox{if } x = y, \\
            0 & \mbox{otherwise.}
            \end{array} \right.
\end{equation}

We denote by $\tau_1 = \tau_1(G)$ the mixing time of the random walk in total variation,
by $\phi = \phi(G)$ the Cheeger constant of the graph $G$
(see subsection \ref{subsec:background} or \cite{Al}).
We will denote by $B_d(n)$ the graph $G = (V,E)$ where
\[
\begin{array}{ll} V = \{-n,\ldots,n\}^d,
& E = \{\left( (x_1,\ldots,x_d) , (y_1,\ldots,y_d) \right) \in V \times V : \sum_{i=1}^d |x_i - y_i| = 1\}.
\end{array}
\]
It is well known that $\tau_1(B_d(n)) = \Theta(n^2)$ and $ \phi(B_d(n)) = \Theta(n^{-1})$.

Let $p_c(\Z^d)$ be the critical parameter for bond percolation in $\Z^d$
(see \cite{g} or subsection \ref{subsec:background}). Fix $p > p_c (\Z^d)$,
and denote by ${\C} = {\C}_d(n) = {\C}_d(n,p)$ the largest cluster inside $B_d(n)$.
Thus ${\C}_d(n)$ is the open component of $B_d(n)$ which has the maximal number of edges.

The following result proves the stability of the mixing time and the
Cheeger constant under percolation:

\begin{Theorem}
\label{thm:Z2}
If $d \geq 2$ and $p > p_c(\Z^d)$, then there exist constants
$0 < c = c(p) < C = C(p)$ such that
\begin{equation} \label{eq:Z2_mix}
\lim_{n \to \infty} \P_p[cn^2 < \tau_1({\C}_d(n)) < C n^2] = 1.
\end{equation}
Similarly, if $d \geq 2$ and $p > p_c(\Z^d)$, then there exist constants
$0 < c = c(p) < C = C(p)$ such that
\begin{equation} \label{eq:Z2_phi}
\lim_{n \to \infty} \P_p[\frac{c}{n} < \phi({\C}_d(n)) < \frac{C}{n}] = 1.
\end{equation} \end{Theorem}


The upper bound on the mixing time is achieved via an estimate of a weighted variant of the Cheeger
constant which was introduced by Lov\`asz and Kannan \cite{lk}.
In order to obtain tight bounds the Lov\`asz-Kannan method,
we study geometrical properties of the percolation cluster, using in particular bootstrap and renormalization.
We conclude this section by recalling some background from percolation theory and
from the theory of finite markov chains.
Section \ref{sec:proofs} contains the proofs of these theorems,
and Section \ref{sec:dis} contains some remarks, conjectures and open problems.

We find it useful to use the following terminology
(analogously to that used in the theory of random graphs).
Let $\G_n$ be an event describing a property of $\C_d(n)$.
We will say that the event holds asymptotically almost surely (abbreviated a.a.s.),
if $\lim_{n \to \infty} \P_p[\G_n] = 1$.

\subsection{Background} \label{subsec:background}

\noindent
{\bf Bernoulli bond percolation.}
In Bernoulli bond percolation on $\Z^d$,
the edges of $\Z^d$ are open (respectively closed)
with probability $p$ (respectively $1-p$) independently.
The corresponding product measure on the configurations of edges
is denoted by $\P_p$ or just $\P$. Let $\C(v)$ be the (open) cluster of $v$.
In other words, $\C(v)$ is the maximal connected component of open edges in
$\Z^d$ containing $v$.

We write
$$
\theta^v(p)=\P_p\big\{\C(v)\mbox{ is infinite} \big\}.
$$
Since $\Z^d$ is transitive, we may write $\theta(p)$ for $\theta^v(p)$.
If $\C(v)$ is infinite for some $v$, we say that {\it percolation} occurs.
We refer the reader to \cite{g} for more background.

A particular property of super-critical percolation that we use below is the following
\begin{Proposition} \label{prop:trivial}
If $p > p_c(\Z^d)$ then there exists $p' = p'(p,d) > 0$ such that a.a.s.
the number of open edges in $\C_d(n)$ is at least $p' d (2n + 1)^d$.
Moreover, there exists a constant $c > 0$ s.t. a.a.s.
$\C_d(n) \cap B_d(n - c \log n)$ is the
intersection of the infinite percolation cluster with $B_d(n - c \log n)$.
\end{Proposition}

\proofs
Let $\theta_e(p)$ be the probability that an edge $e$ belongs to the infinite cluster
(clearly the definition does not depend on the specific edge $e$).
It is clear that since $\theta(p) > 0$ so does $\theta_e(p) > 0$.
Consider $B_d(n)$ as a subgraph of $\Z^d$.
By ergodicity if follows that a.a.s. there are at least
$0.5 \theta_e(p) d (2n + 1)^d$ edges in $B_d(n)$ which belong to the infinite percolation cluster.
By \cite{AP} it follows that there exists some constant $c > 0$ such that a.a.s.
if two edges $e_1,e_2  \in B_d(n - c \log n)$ belong to the infinite percolation cluster,
then they are connected inside $B_d(n)$.

It remains to be shown that in $B_d(n)$ there is at most one connected
component of size larger than $0.5 \theta_e(p) d (2n + 1)^d$.
This follows from the fact (see e.g. \cite{g}) that when $p > p_c$
the probability that a vertex $v$ which is not in the infinite cluster belongs
to a connected cluster of size larger than $k$ is bounded by $\exp(- \alpha k^{(d-1)/d})$ for some $\alpha > 0$.

$\QED$

\medskip
\noindent
{\bf Mixing and relaxation times; Cheeger constant.}
We follow \cite{Al} for basic notations and definitions.
Consider the random walk on the graph $G = (V,E)$ with transition kernel (\ref{eq:transitions})
as a reversible markov chain.
Note that $\pi$, the stationary distribution for the chain satisfies
\begin{equation} \label{eq:pi}
\pi(x) = \frac{\deg(x)}{\sum_{y \in C^d_n \deg(y)}}.
\end{equation}
where $\deg(y)$ is the degree of $y$ in $G$.
Similarly, the probability of an edge $e$, denoted by $Q(e)$, is uniform for all edges of $G$.
For two sets $A$ and $B$ write $Q(A,B) = \sum_{e = (x,y), x \in A, y \in B} Q(e)$.
Let the eigenvalues of $Q$ (\ref{eq:transitions}) be $0 = \lam_1 > \lam_2 \geq \cdots$.
We let the {\em spectral gap} of the random walk on $G$ be $-\lam_2$, and the {\em relaxation time}
of the random walk be $\tau_2 = -\lam_2^{-1}$.

For two distribution measures $\mu$ and $\nu$ on the same discrete space.
The {\em total-variation distance}, $d_V(\mu,\nu)$, between $\mu$ and $\nu$ is defined as
\[
d_V(\mu,\nu) = \frac{1}{2} \sum_{x} |\mu({x}) - \nu({x})|
= \sup_{A} |\mu(A) - \nu(A)|.
\]
Consider again the random walk on $G$.
Denote by $\pi^t_x$ the distribution measure of the walk started at $x$ at time $t$.
The {\em mixing time} of the random walk, $\tau_1$, is defined as
\[
\tau_1 = \inf \{ t : \sup_{x,y} d_V(\pi^t_x,\pi^t_y) \leq e^{-1}\}.
\]

Usually it is harder to estimate $\tau_1$ than it is to estimate $\tau_2$.
However, in general, the following relation holds:
\begin{equation} \label{eq:tau_1_and_tau_2}
\tau_2 \leq \tau_1 \leq \tau_2 \left(1 + \frac{1}{2} \log \frac{1}{\min_x \pi(x)} \right)
\end{equation}
(see e.g. \cite{Al} Lemma 23).

A geometric tool which is used in order to bound relaxation times is Cheeger inequality.
For a set $A$ we define its conductance as
\begin{equation} \label{eq:conductacne}
\phi_A = \frac{Q(A,A^c)}{\pi(A) \pi(A^c)}.
\end{equation}
Let $\phi$ be the Cheeger constant:
\begin{equation} \label{eq:def_phi}
\phi = \inf_A \phi_A = \inf_{A : \pi(A) \leq 1/2} \phi_A.
\end{equation}
Cheeger inequality states that
\begin{equation} \label{eq:cheegerinq}
\tau_2 \leq 8 \phi^{-2}.
\end{equation}
(see e.g. \cite{Al} Theorem 40).

In \cite{lk}, Lov\`asz and Kannan introduced the following variant of the Cheeger constant.
For $0 < x \leq 1/2$, let
\begin{equation} \label{eq:def_phix}
\phi(x) = \min_{\{A : 0 < \pi(A) \leq x\}} \phi_A.
\end{equation}
Then it is shown in \cite{lk} that:
\begin{equation} \label{eq:lk}
\tau_1 \leq 32 \int_{\min_i \pi_i}^{1/2} \frac{1}{x \phi^2(x)}.
\end{equation}

\medskip
\noindent
{\bf Acknowledgement:}
Thanks to Oded Schramm and Prasad Tetali for helpful discussions.

\section{Mixing times} \label{sec:proofs}
\subsection{Lower bounds}

We start by proving the upper bound on the Cheeger constant and the lower bound on the mixing time.
\begin{Lemma} \label{lem:upper}
For all $d \geq 2$ and $p > p_c(\Z^d)$ there exists a constant
$c = c(p,d) > 0$ such that a.a.s. it holds for $\C_d(n)$ that
if $x$ satisfies $\min_v \pi(v) \leq x \leq 1/2$, then
\begin{equation} \label{eq:lowphi}
\phi(x) < \frac{c}{n x^{1/d}}.
\end{equation}
\end{Lemma}

\proofs
By Proposition 1.2. there exits $p'(p)$ such that a.a.s. there are at least
$p' d (2n + 1)^d$ edges belonging to $\C_d(n)$.
In order to prove that (\ref{eq:lowphi}) holds for $\min_v \pi(v) \leq x \leq 1/2$,
it suffices to prove that it holds for $\min_v \pi(v) \leq x \leq q$ for some constant $q$,
as the function $\phi$ is decreasing. Let $r$ be such it is possible to cover at least
$(1 - p'/2)$ of the edges of $B_d(n)$ by disjoint translations of $B_d(k)$ for all $k < n/r + 1$.
Let $k < n/r + 1$. It is clear that at least one of the translations $v + B_d(k)$ satisfies
\[
\frac{2}{p'} \left( \frac{2k + 1}{2n + 1} \right)^d \geq
\pi(B_d(k) + v) \geq \frac{p'}{2} \left( \frac{2k + 1}{2n + 1} \right)^d.
\]
(the inequality on the left is satisfied by all translations).

On the other hand since in $\Z^d$ there are
$2d (2k + 1)^{d-1}$ edges going out of $v + B_d(k)$,
\[
Q(v + B_d(k), \C_d(n) \setminus v + B_d(k)) \leq
\frac{2d (2k+1)^{d-1}}{p' d (2n + 1)^d} = \frac{2 (2k + 1)^{d-1}}{p' (2n + 1)^d}.
\]
It therefore follows that
\[ \phi_{v + B_d(k)} \leq \frac{4}{{p'}^2 (2k + 1)} \leq
 \frac{c'}{n \pi^{1/d}(v + B_d(k))} \leq
 \frac{c}{n} \left(\frac{p'}{2} \left( \frac{2n + 1}{2k + 3} \right)^d \right)^{1/d}
\]
for some constants $c'$ and $c$.
Choosing $k$ to be the maximal such that
\[ x \geq \frac{2}{p'} \left( \frac{2k + 1}{2n + 1} \right)^d \]
we obtain the required result with $q = p' r^{-d}/2$. $\QED$

In order to prove the lower bound on the mixing time, we will use the following lemma.
Let $\pi$ be the stationary distribution for the simple random walk on the graph $G$ (see (\ref{eq:pi})).
Let $f : V \to \R$ be a function. We write $\pi[f] = \sum_{v \in V}
\pi(v) f(v)$ for the expected value
of $f$ with respect to $\pi$.
\begin{Lemma} \label{lem:testdist}
Let $G=(V,E)$ be a finite graph.
For each $v \in V$, let $D_v : V \to R$ be defined by $D_v(x) = D(v,x)$
where $D$ is the graph metric distance. Then
\[
\tau_2 \geq \max_v (\pi[D_v^2] - \pi^2[D_v]).
\]
\end{Lemma}

\proofs
Take $D_v$ as a test function in the extremal characterization of $\tau_2$ (see e.g. \cite{Al}):
\[
\tau_2 = \sup\{\frac{\pi[g^2] - \pi^2[g]}{\dir(g,g)} : \dir(g,g) \neq 0\},
\]
where $\dir$ is the Dirichlet form:
\[
\dir(g,g) = \frac{1}{2} \sum_{u} \sum_{w \neq u}
\pi[u] Q_{u,w} (g(u) - g(w))^2,
\]
and note that if $Q_{v,w} \neq 0$, then
$D_v(u) - D_v(w) \in \{-1,0,1\}$. $\QED$

\begin{Lemma}
For all $d \geq 2$,
there exists a constant $c > 0$ such that $\lim_{n \to \infty} \P_p[\tau_1 \geq c n^2] = 1$.
\end{Lemma}

\proofs
Without loss of generality assume that $0$ belong to the largest open cluster inside $B_d(n)$.
From \cite{AP} it follows that there exist $a > 0, b > 0$
such that a.a.s. if $x , y \in B_d(n - b \log n) \cap \C_d(n)$, then
$D(x,y) \leq a |x - y|_1$.

From proposition \ref{prop:trivial} it follows that a.a.s.
there are at least $c' n^d$ vertices $x$ with
$|x|_1 \leq n/4a$ and at least $c' n^d$ vertices with
$(1 - 1/(8a)) n \geq |x|_1 \geq (1 - 1/(4a))n$ for some constant $c > 0$.
Therefore, there are at least $c' n^d$ vertices $x$ with $D_0(x) \leq n/4$
and at least $c' n^d$ vertices with $D_0(x) \geq 3n/4$.

Applying lemma \ref{lem:testdist} with the function $D_0$ it follows that a.a.s.
$\tau_2 \geq c n^2$ for some constant $c > 0$.
Using the lower bound (\ref{eq:tau_1_and_tau_2}) on $\tau_1$ in terms of $\tau_2$
we achieve the desired conclusion. $\QED$

\subsection{Bootstrap}
In this subsection we show how Theorem \ref{thm:phi} below implies
the upper bound in Theorem \ref{thm:Z2}.
The proof of Theorem \ref{thm:phi} is given in the following subsections.

\begin{Theorem} \label{thm:phi}
For all $d \geq 2$, and all $p > p_c(\Z^d)$,
there exist constants $c_1 = c_1(d,p) > 0$ and $c_2 = c_2(d,p) > 0$
such that $\C_d(n)$ satisfies a.a.s. that for sets $A$ such that
$A$ and $A^c$ are connected and
$\frac{1}{2} \geq \pi(A) \geq \frac{c_1 \log^\frac{d}{d-1} n}{n^d}$,
\begin{equation} \label{eq:phi_x_Z2}
\phi_A \geq \frac{c_2}{n \pi^{1/d}(A)}.
\end{equation} \end{Theorem}

We will also utilize the following lemma which has the same proof as Lemma 36 of \cite{Al}.
\begin{Lemma} \label{comp}
For all $x$, the minimum in (\ref{eq:def_phix}) is obtained for a set $A$
such that $A$ is connected.
The minimum at (\ref{eq:def_phi}) is obtained for a set $A$ such that $A$ and $A^c$ are connected.
\end{Lemma}

For $x$, the set $A$ for which the value $\phi(x)$ is obtained is connected.
However, it may be the case that the complement of the set is not connected.
We bootstrap in the lemma below in order to prove that the estimates in
Theorem \ref{thm:phi} suffice for our purposes.
\begin{Lemma} \label{lem:phi} For all $d \geq 2$, and all $p > p_c(\Z^d)$,
there exist constants $c_1 = c_1(d,p) > 0$ and $c_2 = c_2(d,p) > 0$
such that $\C_d(n)$ satisfies a.a.s. that for all
$\frac{1}{2} \geq x \geq \frac{c_1 \log^\frac{d}{d-1} n}{n^d}$,
\begin{equation} \label{eq:phi_x_Zd}
\phi(x) \geq \frac{c_2}{n x^{1/d}}.
\end{equation}
\end{Lemma}
\proofs
By Lemma \ref{comp} the claim holds for $x=1/2$.
Thus by the monotonicity of $\phi$ in $x$ it follows that for any $q > 0$
by increasing the value of $c_2$, we obtain that a.a.s. (\ref{eq:phi_x_Zd}) holds
for all $q \leq x \leq 1/2$.
We will therefore prove that (\ref{eq:phi_x_Zd}) holds a.a.s. for all
$q \geq x \geq \frac{c_1 \log^\frac{d}{d-1} n}{n^d}$, where $q$ is determined below.
Assume that $A$ is the set for which $\phi(x) = \phi_A$.
By Lemma \ref{comp} the set $A$ is connected. If $\C_d(n) \setminus A$ is connected, we are done,
so we assume the contrary.
Note that by Lemma \ref{lem:upper} a.a.s. $\phi_A \leq \frac{c'}{n x^{1/d}}$
for some constant $c'$. In particular $Q(A,\C_d(n) \setminus A) \leq \frac{c' x^{1 - 1/d}}{n}$.
By the assumption that $\C_d(n) \setminus A$ is not connected we may write
$C_d(n) \setminus A$ as the union of disconnected components $A_1,\ldots,A_r$,
where $r \geq 2$. Let $i$ be the index for which $\pi(A_i)$ is maximized.
Without loss of generality we may assume that $\pi(A_i) \leq 1/2$
(otherwise repeat the argument below for $\C_d(n) \setminus A$).
\begin{Claim} \[
\pi(A_i) \geq \pi(A).
\]
\end{Claim}
\proofs
Assume the contrary and that $q \leq 1/4$.
Note that by Lemma \ref{comp} the set for which the value of $\phi$ is obtained
is connected and has a connected complement.
In particular, for some constant $c''$, a.a.s. $\phi \geq c''/n$.
We may find a sub-collection $I \subset \{1,\ldots,r\}$
such that if $B = A \cup_I A_i$, then $3/4 > \pi(B) > 1/4$.
Note however, that
$Q(B,\C_d(n) \setminus B) = Q(A, \C_d(n) \setminus A)$ and that
by Lemma \ref{lem:upper} a.a.s. $\phi_A \leq \frac{c'}{n x^{1/d}}$ for some constant $c'$.
This implies that
\[ \frac{c''}{n} \leq \phi \leq \phi_B \leq 8 Q(A, \C_d(n) \setminus A) \leq
\frac{8 c'}{n} x^{1 - 1/d} \leq \frac{8 c'}{n} q^{1 - 1/d}.
\]
Thus, when $q$ is sufficiently small, we obtain a contradiction and the proof follows.
$\QED$.

Let $A'$ be the ($\pi$) smallest set among $A_i$ and $C_d(n) \setminus A_i$,
so that $\pi(A') \leq 1/2$.
Note that $A'$ is connected, has a connected complement and satisfies
$1 - \pi(A) \geq \pi(A') \geq \pi(A)$.
It follows by Theorem \ref{thm:phi} that
$\phi_{A'} \geq \frac{c_2}{n \pi^{1/d}(A')}$.
However, this implies that
\[
Q(A,\C_d(n) \setminus A) =
Q(A', \C_d(n) \setminus A') \geq c_2 \frac{\pi^{1 - 1/d}(A')(1 - \pi(A'))}{n}
\geq c_2 \frac{\pi^{1 - 1/d}(A)(1 - \pi(A))}{n}
\]
and we obtain that
\[
\phi_A \geq \frac{c_2}{n x^{1/d}},
\] as needed.
$\QED$

\prooft{of the upper bound in theorem~\ref{thm:Z2}}
We assume that (\ref{eq:phi_x_Zd}) holds.
If the set $A$ satisfies $\pi(A) \leq \frac{c_1 \log^\frac{d}{d-1} n}{n^d}$,
then since $\C_d(n)$ is connected it follows that $Q(A,A^c) \geq 1/(2 d n^d)$
and therefore
\begin{equation} \label{eq:small_sets}
\phi_A \geq \frac{Q(A,A^c)}{\pi(A)} \geq \frac{1}{2 d c_1 \log^\frac{d}{d-1} n}.
\end{equation}
Since by Lemma \ref{comp} the set $A$ which achieves
the minimum at the definition of the Cheeger constant (\ref{eq:def_phi}) is connected,
we obtain by (\ref{eq:small_sets}) and (\ref{eq:phi_x_Z2})
that
\[
\lim_{n \to \infty} \P_p[\phi \geq \frac{2^{1/d} c_2}{n}] = 1.
\]
Moreover, by (\ref{eq:small_sets}) and (\ref{eq:phi_x_Zd}) we obtain that a.a.s.
for all $\frac{1}{2} \geq x \geq \frac{c_1 \log^\frac{d}{d-1} n}{n^d}$,
\[
\phi(x) \geq \min\{\frac{1}{2 d c_1 \log^\frac{d}{d-1} n},\frac{c_2}{n x^{1/d}}\}.
\]

Thus, by (\ref{eq:lk}),
we see that a.a.s.
\[
\tau_1 \leq 32 \int_{\min_v \pi(v)}^{\frac{1}{2}} \frac{1}{t \phi^2(t)} dt \leq
 32 \int_{\frac{1}{2 d n^d}}^{\frac{1}{2}} \frac{4 d^2 c_1^2 \log^\frac{2d}{d-1} n}{t} dt +
 32 \int_{0}^{\frac{1}{2}} \frac{n^2}{c_2^{2} t^{1 - 2/d}} dt \leq C n^2.
\]
for some constant $C = C(d,p)$ as needed. $\QED$

\subsection{Upper bound for $d=2$}

The main tool in the proof will be the following
large deviation result by Kesten \cite{k1}:
\begin{Lemma} \label{lem:large_kesten}
Consider i.i.d. percolation $\{X(e)\}$ on the edges of $\Z^2$ where
\[
\P_p[X(e) = 0] = 1 - \P_p[X(e) = 1] = p < 1/2 \,\,(= p_c(\Z^2)).
\]
For two point $x$ and $y$, let $D(x,y)$ be their distance in the first passage percolation
model:
\[
D(x,y) = \min \{\sum_{i=1}^k X(e_i) : e_1,\ldots,e_k \mbox{ is a path connecting } x \mbox{ to } y \}. \] Then, there exist constants $a > 0$ and $b > 0$ such that for all points $x$ and $y$, \[
\P_p[D(x,y) \leq a |x-y|_1] \leq \exp(- b |x-y|_1),
\]
where $|(x_1,y_1)-(x_2,y_2)|_1 = |x_1 - y_1| + |x_2 - y_2|$.
\end{Lemma}

\prooft{of theorem~\ref{thm:phi}, $d = 2$}

We will use the following dual first-passage percolation model.
Take $\Z^2$ and draw the dual lattice $\Z_{\ast}^2$.
Each edge $e$ of $\Z^2$ crosses a unique edge $e^{\ast}$ of $\Z_{\ast}^2$.
If the edge $e$ is closed, set $X(e^{\ast}) = 0$; otherwise, set $X(e^{\ast}) = 1$.
Note that $\P_p[X(e^{\ast}) = 0] = 1 - p < 1/2$ and we may therefore apply
Lemma \ref{lem:large_kesten} to $\{X(e^{\ast})\}$.
In particular, we obtain that if $c'$ is large, then a.a.s.
for all pairs of points $x^{\ast}$ and $y^{\ast}$ in the dual of $\C_2(n)$ such that
$|x^{\ast}-y^{\ast}|_1 \geq c' \log n$, we have,
\begin{equation} \label{eq:large_kesten}
D(x^{\ast},y^{\ast}) \geq a |x^{\ast} - y^{\ast}|_1.
\end{equation}

Let $A$ be a connected set in $\C_2(n)$ such that
$1/2 \geq t = \pi(A) \geq \frac{c_1 \log^2 n}{n^2}$ and such that
$\C_d(n) \setminus A$ is connected. $B_2(n) \setminus A$ is union of
($\Z^2$) disconnected components $A_1,\ldots,A_r$ where $\C_2(n) \setminus A \subset A_r$.
Let $A' = A \cup \cup_{i=1}^{r-1} A_i$. Then $A'$ and $B_2(n) \setminus A'$ are both connected.
Moreover by Proposition \ref{prop:trivial} a.a.s.
both $A'$ and $B_2(n) \setminus A'$ contain at least $q t n^2$ vertices for some constant $q$.

Let $\gamma$ be the boundary of $A'$ in $B_2(n)$.
In other words, $\gamma$ consists of all the edges $(x,y)$ such that $x \in A', y \notin A'$.
We let $\gamma^{\ast}$ be the path which is obtained by taking the dual edges of the edges of $\gamma$.
Since both $A'$ and $B_2(n) \setminus A'$ are connected, $\gamma^{\ast}$ is connected.

By the isoperimetric inequality for the square in the lattice $\Z^2$ (see \cite{bl}),
$A'$ has $l_1$ diameter at least $c_I \sqrt{q t} n$ for some constant $c_I > 0$.
In other words, there exist two point $x^{\ast}$ and $y^{\ast}$ on $\gamma^{\ast}$
such that $|x^{\ast} - y^{\ast}|_1 \geq c_I \sqrt{q t} n $. Let $|\gamma^{\ast}|$
be the number of edges $e$ on $\gamma^{\ast}$ for which $X(e) = 1$.
Taking $c_1$ sufficiently large it follows by (\ref{eq:large_kesten}) that
a.a.s. $|\gamma^{\ast}| \geq a c_I \sqrt{q t} n$. Since $10 n^2 Q(A,\C_d(n) \setminus A) \geq |\gamma^{\ast}|$,
we obtain that
\[
\phi_A \geq \frac{ a c_I \sqrt{q t}}{10 n t(1-t)} \geq \frac{c_2}{\sqrt{t} n},
\]
for some constant $c_2 = c_2(d,p)> 0$ as needed. $\QED$

\subsection{Upper bound for $d \geq 2$ and $p$ close to $1$}
We now prove the theorem for $d \geq 2$ assuming that $p$ is close to $1$.
For two vertices $v$ and $w$ in $B_d(n)$,
a {\sl cutset} separating $v$ from $w$ is a set $B$ of edges of $B_d(n)$
such that any path in $B_d(n)$ which connects $v$ to $w$ intersects $B$ at least in one edge.
A {\sl minimal} cutset separating $v$ from $w$,
is a cutset which has no proper subset which is also a cutset separating $v$ from $w$.

The following fact is probably well known.
We refer the reader to Babson and Benjamini \cite{bb} for a proof (in a more general setting).
\begin{Lemma} \label{lem:cut_sets}
For $v,w \in B_d(n)$ the number of minimal cutsets of size $m$
separating $v$ from $w$ is bounded by $c(d)^m$ for some constant $c(d)$.
\end{Lemma}

\prooft{ of Theorem \ref{thm:phi}, $d \geq 2$ and $p$ close to $1$}
By Proposition \ref{prop:trivial} a.a.s. the
$\pi$ measure of a subset $A \subset \C_d(n)$
is up to constant the same as the number of vertices in the set divided by
$(2n + 1)^d$.
Therefore, in order to show that (\ref{eq:phi_x_Z2}) holds,
it suffices to show that there exist constants $\hat{c}> 0$ and $\hat{c}_1 > 0$
such that a.a.s. all connected sets $A \subset \C_d(n)$ such that
$\C_d(n) \setminus A$
is also connected and such that the
size of $|A|$ is at least ${\hat{c}}_1 \log^\frac{d}{d-1} n$
and at most $|\C_d(n)|/2$, satisfy that the number of open edges going from
$A$ to $\C_d(n) \setminus A$ is at least $\hat{c} |A|^{\frac{d-1}{d}}$.

For such a set $A$, the set $B_d(n) \setminus A$
is a union of disconnected components
$A_1,\ldots,A_r$ where $\C_d(n) \setminus A \subset A_r$.
Fix a point $v \in A$ and $w \in A_r$ and look at the set $B$ of edges connecting
$A$ to $A_r$ ($=$ the set of edges connecting $A \cup \cup_{i=1}^{r-1} A_i$ to $A_r$).
This is the minimal cutset of edges separating $v$ from $w$.
Moreover, since $|A_r| \geq |A|$, it follows by the isoparametric inequality \cite{bl}
that $|B| \geq c_I |A|^{\frac{d-1}{d}} \geq c_I {\hat{c}_1}^\frac{d-1}{d} \log n$ for some constant $c_I$.

It follows that in order to prove the theorem,
it suffices to show that there exist constants $c_1 > 0, c_2 > 0$ such that
for all minimal cut-sets $|B|$ of size at least $c_1 \log n$,
the number of open edges in $|B|$ is at least $c_2 |B|$.
We denote by $\tilde{B}$ the subset of open edges of set $B$.

We will use a first moment argument.
Applying large deviation estimates, we see that if $c_2 > 0$ is sufficiently small,
and $p < 1$ is sufficiently large, then $\P_p[|\tilde{B}| \leq c_2 |B|] \leq (c(d)+1)^{-|B|}$.

Summing up, and using Lemma \ref{lem:cut_sets}
we see that the probability that there exists any cut-set
$B$ with $|\tilde{B}| \leq c_2 |B|$ and $|B| \geq c_1 \log n$ is bounded by
\[
(2n+1)^{2d} \sum_{s \geq c_1 \log n} c(d)^s (c(d)+1)^{-s} = o(1),
\]
provided that $c_1$ is sufficiently large.
$\QED$

\subsection{Upper bound for $d \geq 2$ and $p > p_c(\Z^d)$} 

The proof uses renormalization and the result for $p$ close to $1$.
The renormalization will produce site percolation with high density of "good" sites.
We will need the following fact
\begin{Proposition} \label{prop:connected}
There exists constants $p^{\ast} < 1,a > 0$ and $c > 0$ such that a.a.s.
for site percolation with parameter $p > p^{\ast}$ all connected sets $A$ in $B_d(n)$
of size at least $c \log n$ have at least $a |A|$ open sites.
\end{Proposition}

\proofs
It is well known that the number of connected sets of size $m$
containing a specified vertex $v$ is bounded by $c(d)^m$ for some constant $c(d)$.
If $p^{\ast} < 1$ is sufficiently large, and $a > 0$ is sufficiently small,
then for each set $A$, the number of open sites in $A$ denoted $|\tilde{A}|$ satisfies
\[
\P_p[|\tilde{A}| \leq a |A|] \leq (c(d)+1)^{-|A|}.
\]
Summing over all sets,
we see that the probability that there exists any connected set of size greater than
$c \log n$ for which $|\tilde{A}| < a |A|$ is bounded by
\[
(2n+1)^d \sum_{s \geq c \log n} c(d)^s (c(d)+1)^{-s} = o(1),
\]
provided that $c$ is sufficiently large.
$\QED$

We have proved Theorem \ref{thm:phi} and therefore Lemma \ref{lem:phi}
for large $p < 1$.
An analogous proof implies the analogous result for site percolation for large $p < 1$.
For convenience we state this result below:
\begin{Lemma} \label{lem:site}
For all $d \geq 2$, there exists $p^{\ast} < 1$ such that for $p > p^{\ast}$,
there exist constants $c_1 = c_1(d,p) > 0$ and $c_2 = c_2(d,p) > 0$ such that a.a.s.
site percolation with parameter $p$ on $B_d(n)$ satisfies for all
$\frac{1}{2} \geq x \geq \frac{c_1 \log^\frac{d}{d-1} n}{n^d}$ that
\begin{equation} \label{eq:phi_site}
\phi(x) \geq \frac{c_2}{n x^{1/d}}.
\end{equation}
\end{Lemma}

For $v \in \Z^d$, we let
\[
Q_N(v) = v + B_d(N) = \{w : |w - v|_{\infty} \leq N\}.
\]
We will slightly abuse the notation by writing $Q_N(v)$
for the induced subgraph on $Q_N(v)$.
We call $v \in (N \Z)^d$ a {\sl good} vertex if the following conditions
hold:
\begin{itemize}
\item
There exists an open cluster which intersects all $d-1$ dimensional faces of the box $Q_{5N/4}(v)$.
\item All connected components of diameter more than $N/10$ inside the box $Q_{5N/4}(v)$ intersect the above cluster.
\end{itemize}
By standard renormalization results (see Proposition 2.1 in Antal and Pisztora \cite{AP})
it follows that for any $p>p_c(\Z^d)$,
the set of good vertices stochastically dominates site percolation with parameter
$p^\ast(N)$ on $(N \Z)^d$, with $\lim_{N\rightarrow\infty}p^\ast(N)=1$.

We take a connected set $A$ in $\C_d(n)$ such that both $A$ and $\C_d(n) \setminus A$ are connected and
such that $C \log^{\frac{d}{d-1}} n \leq |A| \leq |\C_d(n) \setminus A|$.
We will show
that a.a.s. the number of open edges between $A$ and $\C_d(n)$ is at least
$c |A|^{\frac{d-1}{d}}$, where $c$ and $C$ are positive constants to be determined later.
We will thus obtain the required result.

We let $A' = \{v \in (N \Z)^d : |Q_{5N/4}(v) \cap A| \geq N/10\}$.
Since $A$ is connected, $A'$ is a connected set in $(N \Z)^d$.
Moreover, it is clear that $|A|/(2N)^d \leq |A'| \leq |A|$.
Let $A_g$ be the set of good sites in $A'$. By Proposition \ref{prop:connected}
when $p^{\ast}(N)$ is sufficiently large, a.a.s. $|A_g| \geq a |A'| \geq \frac{a}{(2N)^d} |A|$.

It now follows by Lemma \ref{lem:site} that if $C$ is sufficiently large,
then a.a.s. there are at least $c'' |A_g|^{\frac{d-1}{d}} \geq c'
|A|^{\frac{d-1}{d}}$
pairs of good neighbors $u$ and $w$ such that $u \in A_g$ and $w \notin A_g$, where $c' > 0$.

We note that each such pair defines an open edge going from $A \cap Q_N(u)$ to
$(\C_d(n) \setminus A) \cap Q_N(w)$.
Moreover each such edge is defined by at most $4 d^2$ pairs $(u,w)$.
It now follows that a.a.s. the edge boundary of $A$ is of size at least
$\frac{c'}{4 d^2} |A|^{\frac{d-1}{d}}$ as needed.
$\QED$

\section{Further remarks} \label{sec:dis}
\subsection{Coupling}
For simple random graph models it is easy to bound from above the mixing time by constructing explicit coupling.
We give two examples below

\begin{enumerate}

\item
Let $G(n,p)$ be the random graph model for fixed $p$.
It is easy to see that every two vertices have at least $\frac{p^2 n}{2}$
joint neighbors with probability going to $1$ as $n \to \infty$.
Therefore, by coupling we see that $\sup_{x,y} D_V(\P^t_x,\P^t_y) < (1 - p^2)^t$.
So the mixing time is $\Theta(1)$.

\item \label{ex:dembo}
Consider the following perturbation of $B_2(n)$.
To each vertex $v$ of the square attach a pipe of length
$X_v$ where $X_v$ are independent random variables taking the values $0,1$.
As noted to us by Amir Dembo, one can use the usual reflecting coupling on the square $B_2(n)$,
in order to show that the mixing time is $\Theta(n^2)$.
Indeed, let $x$ and $y$ be two vertices and consider random walks starting at $x$ and $y$.
Always delay one of the walks in order that the two walks make steps in $B_2(n)$ simultaneously.
Whenever the two walks make steps in $B_2(n)$, use the usual reflecting coupling.
This example can be generalized to any dimension and the assumption on $X_v$ may be replaced
by a weaker moment assumption. In a previous draft of this paper,
we had a more complicated result in the same spirit.

\end{enumerate}

We think it is an interesting challenge to try a variant of the last argument in order
to show that the mixing time on $C^d_n$ is $O(n^2)$ for all $d$.
One approach of implicit construction of such coupling is to use some kind
of central limit theorem in order to bring the walks closer and closer. Unfortunately
the present form of the CLT on super critical percolation cluster (De Masi et. al. (1989))
provides no estimates on the convergence rate and therefore no bound on the coupling time of $\C_d(n)$.

\subsection{Final comments}
\begin{enumerate}
\item
In this note $ p > p_c(\Z^d)$, is fixed, what is the dependence on $p$ of the mixing time?
What is the mixing times on the critical cluster?

\item
One can consider the mixing time of random walks on percolation clusters on other graphs.
With Nick Wormald (in preparation), it is shown that the mixing time for simple random walk on
$G(n, c/n), c >1 $ is poly-logarithmic in $n$.

\item
It is natural to ask if the cover time of simple random walk is also robust under perturbations.
Using the "stretched lattice" representation of the percolation cluster, we show (in preparation) that the cover
time of $\C_d(n)$ for $d \geq 3$ is $\Theta( n^d \log^2 n)$ compared with $\Theta(n^d \log n)$ for
the cover time of $B_d(n)$.
\end{enumerate}


\vspace{.2in}
\noindent {\sc Itai Benjamini} \newline
The Weizmann Institute  and Microsoft Research \newline
{\tt itai@wisdom.weizmann.ac.il}\newline


\vspace{.1in}
\noindent {\sc Elchanan Mossel} \newline
Microsoft research \newline
{\tt mossel@microsoft.com}\newline

\vspace{.25in}
\noindent


\begin{thebibliography}{19}




\bibitem{Al}
D. Aldous and J. A. Fill, Reversible Markov chains and random walks on graphs,
{\em book in preperation}. (2000)

\bibitem{AP} Antal, P. and Pisztora, A. (1996)
On the chemical distance in supercritical Bernoulli percolation,
{\sl Ann. Probab.} {\bf 24}, 1036--1048.


\bibitem{bb}
E. Babson, and I. Benjamini, Cut sets and normed cohomology with applications to percolation.
{\it Proc. Amer. Math. Soc.} {\bf 127} (1999), no. 2, 589--597.

\bibitem{bl}
Bollob\'{a}s, B, and Leader, Imre, Edge-isoperimetric inequalities in the grid.
{\it Combinatorica} {\bf 11} (1991), no. 4, 299--314.

\bibitem{cls}
T. Liggett, R. Schonmann, R. and A. Stacey,
Domination by product measures.
{\it Ann. Probab.} {\bf 25} (1997), no. 1, 71--95.

\bibitem{dfgw}
A. De Masi, P. Ferrari, S.Goldstein, W. Wick,
An invariance principle for reversible Markov processes.
Applications to random motions in random environments.
{\it J. Statist. Phys.} 55 (1989), no. 3-4, 787--855.

\bibitem{g}
G. Grimmett,  Percolation. Second edition.
Springer-Verlag, Berlin, 1999. xiv+444


\bibitem{k1}
H. Kesten, On the time constant and path length of first-passage percolation.
{\it Adv. in Appl. Probab.} 12 (1980), no. 4, 848--863.


\bibitem{lk}
L. Lov\`asz and R. Kannan,
Faster mixing via average conductance {\it Proc. 1995 ACM STOC}.


\end{thebibliography}
\end{document}